 \documentclass[12 pt]{amsart}
\usepackage{comment}
\usepackage{enumerate}
\usepackage{amssymb, amsmath}

\usepackage{url}

\DeclareMathOperator{\CM}{\mathrm CM}

\DeclareMathOperator{\norm}{\mathrm norm}
\DeclareMathOperator{\Norm}{\mathrm Norm}

\DeclareMathOperator{\rank}{\mathrm rank}
\DeclareMathOperator{\Reg}{\mathrm Reg}

\DeclareMathOperator{\Tor}{\mathrm Tor}

\DeclareMathOperator{\Vol}{\mathrm Vol}

\date{}

\begin{document}

 \newtheorem{theorem}{Theorem}[section]
 \newtheorem{lemma}{Lemma}[section]
 \newtheorem{corollary}{Corollary}[section]
 \newtheorem{conjecture}{Conjecture}
 \newtheorem{proposition}{Proposition}[section]
 \newtheorem{definition}{Definition}
 \newcommand{\mc}{\mathcal}
 \newcommand{\A}{\mc A}
 \newcommand{\B}{\mc B}
 \newcommand{\eC}{\mc C}
 \newcommand{\D}{\mc D}
 \newcommand{\E}{\mc E}
 \newcommand{\F}{\mc F }
 \newcommand{\G}{\mc G}
 \newcommand{\hH}{\mc H}
 \newcommand{\I}{\mc I}
 \newcommand{\J}{\mc J}
 \newcommand{\K}{\mc K}
 \newcommand{\eL}{\mc L}
 \newcommand{\M}{\mc M}
 \newcommand{\eN}{\mc N}
 \newcommand{\pP}{\mc P}
 \newcommand{\qq}{\mc Q}
 \newcommand{\eS}{\mc S}
 \newcommand{\eU}{\mc U}
 \newcommand{\V}{\mc V}
 \newcommand{\X}{\mc X}
 \newcommand{\Y}{\mc Y}
 \newcommand{\C}{\mathbb{C}}
 \newcommand{\R}{\mathbb{R}}
 \newcommand{\N}{\mathbb{N}}
 \newcommand{\bP}{\mathbb P}
 \newcommand{\Q}{\mathbb{Q}}
 \newcommand{\T}{\mathbb{T}}
 \newcommand{\Z}{\mathbb{Z}}
 \newcommand{\fA}{\mathfrak A}
 \newcommand{\fB}{\mathfrak B}
 \newcommand{\fD}{\mathfrak D}
 \newcommand{\fd}{\mathfrak d}
 \newcommand{\fE}{\mathfrak E}
 \newcommand{\ff}{\mathfrak F}
 \newcommand{\fI}{\mathfrak I}
 \newcommand{\ffi}{\mathfrak i}
 \newcommand{\fJ}{\mathfrak J}
 \newcommand{\fj}{\mathfrak j}
 \newcommand{\fP}{\mathfrak P}
 \newcommand{\fp}{\mathfrak p}
 \newcommand{\fQ}{\mathfrak Q}
 \newcommand{\fq}{\mathfrak q}
 \newcommand{\fS}{\mathfrak S}
 \newcommand{\fU}{\mathfrak U}
 \newcommand{\fu}{\mathfrak u}
 \newcommand{\fV}{\mathfrak V}
 \newcommand{\fv}{\mathfrak v}
 \newcommand{\fb}{f_{\beta}}
 \newcommand{\fg}{f_{\gamma}}
 \newcommand{\gb}{g_{\beta}}
 \newcommand{\vphi}{\varphi}
 \newcommand{\vep}{\varepsilon}
 \newcommand{\bo}{\boldsymbol 0}
 \newcommand{\bone}{\boldsymbol 1}
 \newcommand{\ba}{\boldsymbol a}
 \newcommand{\bb}{\boldsymbol b}
 \newcommand{\bc}{\boldsymbol c}
 \newcommand{\be}{\boldsymbol e}
 \newcommand{\bk}{\boldsymbol k}
 \newcommand{\bell}{\boldsymbol \ell}
 \newcommand{\bm}{\boldsymbol m}
 \newcommand{\bn}{\boldsymbol n}
 \newcommand{\balpha}{\boldsymbol \alpha}
 \newcommand{\bgamma}{\boldsymbol \gamma}
 \newcommand{\bt}{\boldsymbol t}
 \newcommand{\bu}{\boldsymbol u}
 \newcommand{\bv}{\boldsymbol v}
 \newcommand{\bx}{\boldsymbol x}
 \newcommand{\bX}{\boldsymbol X}
 \newcommand{\bwy}{\boldsymbol y}
 \newcommand{\Bbeta}{\boldsymbol \beta}
 \newcommand{\bxi}{\boldsymbol \xi}
 \newcommand{\bbeta}{\boldsymbol \eta}
 \newcommand{\bw}{\boldsymbol w}
 \newcommand{\bz}{\boldsymbol z}
 \newcommand{\bzeta}{\boldsymbol \zeta}
 \newcommand{\hmu}{\widehat \mu}
 \newcommand{\oK}{\overline{K}}
 \newcommand{\oKt}{\overline{K}^{\times}}
 \newcommand{\oQ}{\overline{\Q}}
 \newcommand{\oq}{\oQ^{\times}}
 \newcommand{\oQt}{\oQ^{\times}/\Tor\bigl(\oQ^{\times}\bigr)}
 \newcommand{\ot}{\Tor\bigl(\oQ^{\times}\bigr)}
 \newcommand{\h}{\frac12}
 \newcommand{\hh}{\tfrac12}
 \newcommand{\dx}{\text{\rm d}x}
 \newcommand{\dbx}{\text{\rm d}\bx}
 \newcommand{\dy}{\text{\rm d}y}
 \newcommand{\dby}{\text{\rm d}\bwy}
 \newcommand{\dmu}{\text{\rm d}\mu}
 \newcommand{\dnu}{\text{\rm d}\nu}
 \newcommand{\ds}{\text{\rm d}s}
 \newcommand{\trho}{\widetilde{\rho}}
 \newcommand{\dtrho}{\text{\rm d}\widetilde{\rho}}
 \newcommand{\drho}{\text{\rm d}\rho}
 \newcommand{\dxi}{\text{\rm d}\xi}

\title[Heights]{Independent relative units of low height}
%\date{\vspace{-5ex}}
\author{Shabnam~Akhtari}
\author{Jeffrey~D.~Vaaler}
%\thanks{}
\date{}
\address{Department of Mathematics, University of Oregon, Eugene, Oregon 97403 USA}
\email{akhtari@uoregon.edu}
\address{Department of Mathematics, University of Texas, Austin, Texas 78712 USA}
\email{vaaler@math.utexas.edu}
\allowdisplaybreaks
\numberwithin{equation}{section}

%%%%%%%%%%%%%%%%%%%%%%%%%%%%%%%%%%%%%%%%%%%%%%%%%%%%%%%%%%%%%%
\begin{abstract}  
We  prove inequalities that compare the relative regulator of an extension of number fields with a product of heights of 
multiplicatively independent relative units.
\end{abstract}

\maketitle

%%%%%%%%%%%%%%%%%%%%%%%%%%%%%%%%%%%%%%%%%%%%%%%%%%%%%%%%%%%%%%%
\section{Introduction}\label{Intro}

In \cite[Theorem 1.2]{akhtari2016} we proved the existence of a maximal collection of independent $S$-units 
such that the product of their heights was bounded by a multiple of the $S$-regulator.  Here we establish an analogous result
that bounds the product of the heights of a maximal collection of independent relative units by a multiple of the relative regulator.

We assume that $k$ and $l$ are algebraic number fields such that
\begin{equation}\label{short277}
\Q \subseteq k \subseteq l \subseteq \oQ,
\end{equation}
where $\oQ$ is a fixed algebraic closure of $\Q$.  We write $r(k)$ for the rank of the unit group $O_k^{\times}$, 
and $r(l)$ for the rank of the unit group $O_l^{\times}$.  To avoid degenerate situations we assume that
\begin{equation}\label{short298}
1 \le r(k) < r(l).
\end{equation}
The inequality on the left of (\ref{short298}) implies that $k$ is not $\Q$, and $k$ is not an imaginary, quadratic extension of $\Q$.
The inequality on the right of (\ref{short298}) implies that either $l$ is not a $\CM$-field, or if $l$ is a $\CM$-field then $k$ is not the 
unique, maximal totally real subfield of $l$ (see \cite[Corollary 1 of Proposition 3.20]{narkiewicz2010}).
As roots of unity do not play a significant role in our results, it will be convenient to work in the quotient groups
\begin{equation*}\label{short303}
F_k = O_k^{\times}/\Tor\bigl(O_k^{\times}\bigr),\quad\text{and}\quad F_l = O_l^{\times}/\Tor\bigl(O_l^{\times}\bigr).
\end{equation*}
Then $F_k$ is a free group of rank $r(k)$, and $F_l$ is a free group of rank $r(l)$.  If $\alpha$ belongs to $O_k^{\times}$ and 
$\alpha \Tor\bigl(O_k^{\times}\bigr)$ is the corresponding coset in $F_k$, we have the injective homomorphism
\begin{equation}\label{extra252}
\alpha \Tor\bigl(O_k^{\times}\bigr) \mapsto \alpha \Tor\bigl(O_l^{\times}\bigr).
\end{equation}
The map (\ref{extra252}) allows us to identify $F_k$ with its image in $F_l$, and so we regard $F_k$ as a subgroup of $F_l$.

The restricted norm homomorphisms
\begin{equation*}\label{extra255}
\Norm_{l/k} : O_l^{\times} \rightarrow O_k^{\times},\quad\text{and},
	\quad \Norm_{l/k} : \Tor\bigl(O_l^{\times}\bigr) \rightarrow \Tor\bigl(O_k^{\times}\bigr),
\end{equation*}
induce a homomorphism which we write as
\begin{equation}\label{extra257}
\norm_{l/k} : F_l \rightarrow F_k.
\end{equation}
Then using (\ref{extra257}) we define the subgroup of {\it relative units} in $F_l$ to be the kernel
\begin{equation}\label{extra259}
E_{l/k} = \big\{\alpha \in F_l : \norm_{l/k}(\alpha) = 1\big\}.
\end{equation}
It can be shown (as in \cite[section 3]{akhtari2016}) that $E_{l/k} \subseteq F_l$ is a subgroup such that
\begin{equation*}\label{extra261}
1 \le \rank E_{l/k} = r(l/k) = r(l) - r(k) < r(l).
\end{equation*}
The elements of the group $E_{l/k}$ are {\it relative units}.  We write $\Reg\bigl(E_{l/k}\bigr)$ for the {\it relative regulator} of $l/k$, 
and define this explicitly in (\ref{extra294}) and (\ref{extra301}) (see also \cite[equation (3.5)]{akhtari2016} or 
\cite[equation (1.3)]{costa1991}).  If $\alpha \not= 1$ belongs to $F_k$ then
\begin{equation*}\label{extra265}
\norm_{l/k}(\alpha) = \alpha^{[l:k]} \not= 1,
\end{equation*}
and if $\alpha$ belongs to $E_{l/k}$ then
\begin{equation*}\label{extra266}
\norm_{l/k}(\alpha) = 1.
\end{equation*}
Thus we have
\begin{equation}\label{extra270}
F_k \cap E_{l/k} = \{1\},
\end{equation}
and the image
\begin{equation}\label{extra277}
I_{l/k} = \big\{\norm_{l/k}(\alpha) : \alpha \in F_l\big\} \subseteq F_k
\end{equation}
is a subgroup of full rank $r(k)$.

If $\alpha \not= 0$ is an algebraic number we write $h(\alpha)$ for the Weil height of $\alpha$, and we define this
height in (\ref{B6}).  If $\zeta$ is a root of unity then it is well known that $h(\alpha \zeta) = h(\alpha)$.  Therefore 
the height is well defined as a map
\begin{equation}\label{extra303}
h : \oQt \rightarrow [0, \infty).
\end{equation}
In particular, the height is well defined on cosets of the quotient groups $F_k$ and $F_l$.

Let $\fE \subseteq E_{l/k}$ be a subgroup of maximal rank.  In \cite[Theorem 3.1]{akhtari2016} we proved that each collection of 
multiplicatively independent relative units $\gamma_1, \gamma_2, \dots , \gamma_{r(l/k)}$ contained in $\fE$ satisfies the inequality
\begin{equation}\label{short329}
\Reg\bigl(E_{l/k}\bigr) [E_{l/k} : \fE] \le \prod_{j = 1}^{r(l/k)} \bigl([l : \Q] h(\gamma_j)\bigr).
\end{equation}
Here we show that the inequality (\ref{short329}) is sharp up to a constant that depends only on the rank $r(l/k)$.

\begin{theorem}\label{thmintro1}  Let $k \subseteq l$ be algebraic number fields such that the group $E_{l/k}$ has positive rank 
$r(l/k) = r(l) - r(k)$.  Let $\fE \subseteq E_{l/k}$ be a subgroup of rank $r(l/k)$.  Then there exist multiplicatively independent
relative units $\psi_1, \psi_2, \dots , \psi_{r(l/k)}$ in $\fE$ such that
\begin{equation}\label{short322}
\prod_{j = 1}^{r(l/k)} \bigl([l : \Q] h(\psi_j)\bigr) \le r(l/k)! \Reg\bigl(E_{l/k}\bigr) [E_{l/k} : \fE].
\end{equation}
\end{theorem}

We have remarked that the subgroup $I_{l/k} \subseteq F_k$ defined in (\ref{extra277}) has maximal rank $r(k)$.  Therefore
it follows from \cite[Theorem 1.2]{akhtari2016} that there exist multiplicatively independent elements 
$\beta_1, \beta_2, \dots , \beta_{r(k)}$ in $I_{l/k}$ such that
\begin{equation}\label{short440}
\prod_{i = 1}^{r(k)} \bigl([k : \Q] h(\beta_i)\bigr) \le r(k)! \Reg(k) [F_k : I_{l/k}].
\end{equation}
And it follows from (\ref{short322}) that there exist multiplicatively independent elements $\psi_1, \psi_2, \dots , \psi_{r(l/k)}$
in the group $E_{l/k}$ of relative units such that
\begin{equation}\label{short445}
\prod_{j = 1}^{r(l/k)} \bigl([l : \Q] h(\psi_j)\bigr) \le r(l/k)! \Reg\bigl(E_{l/k}\bigr).
\end{equation}
Thus the units $\beta_1, \beta_2, \dots , \beta_{r(k)}$ are in the image of the homomorphism $\norm_{l/k}$ and the units
$\psi_1, \psi_2, \dots , \psi_{r(l/k)}$ are in the kernel of $\norm_{l/k}$.  In view of (\ref{extra270}) these two sets of multiplicatively 
independent units can be combined.

\begin{corollary}\label{corintro1}  Let $\beta_1, \beta_2, \dots , \beta_{r(k)}$ be multiplicatively independent units in $I_{l/k}$ that 
satisfy {\rm (\ref{short440})}, and let $\psi_1, \psi_2, \dots , \psi_{r(l/k)}$ be multiplicatively independent units in $E_{l/k}$ that satisfy
{\rm (\ref{short445})}.  Then the elements in the set
\begin{equation}\label{short450}
\big\{\beta_1, \beta_2, \dots , \beta_{r(k)}\big\} \cup \big\{\psi_1, \psi_2, \dots , \psi_{r(l/k)}\big\}
\end{equation}
are multiplicatively independent units in the subgroup
\begin{equation*}\label{short452}
I_{l/k} \oplus E_{l/k} \subseteq F_l,
\end{equation*}
and satisfy the inequality
\begin{equation}\label{short455}
\prod_{i = 1}^{r(k)} \bigl([k : \Q] h(\beta_i)\bigr) \prod_{j = 1}^{r(l/k)} \bigl([l : \Q] h(\psi_j)\bigr) \le r(k)! r(l/k)! \Reg(l).
\end{equation}
\end{corollary}

This manuscript is organized as follows.  In section \ref{GeoNum} we prove a general identity that connects the volumes of certain star-bodies
in Euclidean spaces.  In section \ref{heights} we define a generalization of the Schinzel norm that was used in \cite{akhtari2016}, and we record
the identity (\ref{V63}) for the volume of the unit ball attached to the generalized Schinzel norm.  The reason for using the generalized
Schinzel norm when working with relative units can be seen in the basic identity (\ref{short305}) proved in Lemma \ref{lemheight1}.  
Section \ref{proofs} contains our proof of Theorem \ref{thmintro1}, which is an application of Minkowski's theorem on successive minima.
 
%%%%%%%%%%%%%%%%%%%%%%%%%%%%%%%%%%%%%%%%%%%%%%%%%%%%%%%%%%%%%%%
\section{A volume identity}\label{GeoNum}

Let $N$ be a positive integer, and let 
\begin{equation*}\label{U1}
\Delta : \R^N \rightarrow [0, \infty)
\end{equation*}
be a distance function on $\R^N$ in the sense of \cite[Chapter IV]{cassels1971}.  That is, 
\begin{itemize}
\item[(i)]  $\bx \mapsto \Delta(\bx)$ is continuous, and
\item[(ii)]  $\Delta(\xi \bx) = \xi \Delta(\bx)$ for all $0 < \xi < \infty$ and all $\bx$ in $\R^N$.
\end{itemize}
It follows from these conditions that
\begin{equation}\label{U3}
\fS = \big\{\bx \in \R^N : \Delta(\bx) \le 1\big\}
\end{equation}
is a closed {\it star-body} in $\R^N$ (see \cite[Chapter IV, section 1]{cassels1971}).  As $\bo$ is an element of the open star-body
\begin{equation*}\label{U6}
\big\{\bx \in \R^N : \Delta(\bx) < 1\big\} \subseteq \fS,
\end{equation*}
it is obvious that $\Vol_N(\fS)$ is greater than zero.  In this section we will be interested in distance functions $\Delta$ such that
\begin{equation}\label{U9}
\Vol_N(\fS) < \infty.
\end{equation}
If
\begin{equation*}\label{U13}
\chi_{\fS} : \R^N \rightarrow \{0, 1\}
\end{equation*}
is the characteristic function of the closed star-body $\fS$, then (\ref{U9}) implies that
\begin{equation}\label{U16}
\begin{split}
\Vol_N\bigl(\{\bx \in \R^N : \Delta_N(\bx) \le \xi\}\bigr) &= \int_{\R^N} \chi_{\fS}\bigl(\xi^{-1} \bx\bigr)\ \dbx\\ 
										&= \Vol_N(\fS) \xi^N
\end{split}
\end{equation}
for $0 <\xi < \infty$.

\begin{lemma}\label{lemU1}  Let $\Delta$ be a distance function on $\R^N$ and let $\fS$ be the associated star-body defined by
{\rm (\ref{U3})}.  If
\begin{equation}\label{U237}
\Vol_N(\fS) < \infty,
\end{equation}
then at each point $s = \sigma + it$ in the open complex half plane
\begin{equation}\label{U243}
\{s = \sigma + it \in \C : 0 < \sigma\}
\end{equation}  
we have
\begin{equation}\label{U245}
N! \Vol_N(\fS) s^{-N} = \int_{\R^N} \exp\bigl(-s \Delta(\bx)\bigr)\ \dbx.
\end{equation}
\end{lemma}

\begin{proof}  We write
\begin{equation*}\label{U353}
\chi_{\Delta} : \R^N \times (0, \infty) \rightarrow \{0, 1\}
\end{equation*}
for the characteristic function
\begin{equation*}\label{U355}
\chi_{\Delta}(\bx, \xi) = \begin{cases}      1&     \text{if $\Delta(\bx) \le \xi$,}\\
						         0&     \text{if $\xi < \Delta(\bx)$.}\end{cases}
\end{equation*}
If (\ref{U237}) holds and $0 < \sigma$ then
\begin{equation*}\label{U361}
\begin{split}
\int_{\R^N} \int_0^{\infty} \chi_{\Delta}(\bx, \xi) \exp(-\sigma \xi)\ \dxi\ \dbx &= \int_{\R^N} \int_{\Delta(\bx)}^{\infty} \exp(-\sigma\xi)\ \dxi\ \dbx\\
	&= \sigma^{-1} \int_{\R^N} \exp\bigl(-\sigma \Delta(\bx)\bigr)\ \dbx.
\end{split}
\end{equation*}
Applying Fubini's theorem we also get
\begin{equation*}\label{U368}
\begin{split}
\int_{\R^N} \int_0^{\infty} \chi_{\Delta}(\bx, \xi) \exp(-\sigma \xi)\ \dxi\ \dbx 
         &= \int_0^{\infty} \int_{\R^N} \chi_{\Delta}(\bx, \xi)\ \dbx \exp(-\sigma \xi)\ \dxi\\
	&= \Vol_N(\fS) \int_0^{\infty} \xi^N \exp(-\sigma \xi)\ \dxi\\
	&= N! \Vol_N(\fS) \sigma^{- N - 1}.
\end{split}
\end{equation*}
By combining the previous equations we obtain the proposed identity
\begin{equation}\label{U375}
N! \Vol_N(\fS) \sigma^{-N} = \int_{\R^N} \exp\bigl(-\sigma \Delta(\bx)\bigr)\ \dbx
\end{equation}
at each point $\sigma$ on the positive real axis.  

For each complex number $s = \sigma +it$ we have
\begin{equation*}\label{U378}
\biggl|\int_{\R^N} \exp\bigl(- s \Delta(\bx)\bigr)\ \dbx\biggr| \le \int_{\R^N} \exp\bigl(- \sigma \Delta(\bx)\bigr)\ \dbx.
\end{equation*}
Then it follows from (\ref{U375}) that 
\begin{equation}\label{U382}
s \mapsto \int_{\R^N} \exp\bigl(- s \Delta(\bx)\bigr)\ \dbx
\end{equation}
is uniformly bounded on each compact subset of the open half plane (\ref{U243}).  Let $T$ be a closed triangle in (\ref{U243})
and let $\partial T$ be a closed, positively oriented path along the boundary of $T$.  Then it follows from Fubini's theorem that
\begin{equation*}\label{U385}
\begin{split}
\int_{\partial T} \int_{\R^N} &\exp\bigl(- s \Delta(\bx)\bigr)\ \dbx\ \ds\\ 
	                &= \int_{\R^N} \int_{\partial T} \exp\bigl(- s \Delta(\bx)\bigr)\ \ds\ \dbx = 0.
\end{split}
\end{equation*}  
Therefore applying Morera's theorem (see \cite[Theorem 10.17]{rudin1987}) we find that (\ref{U382}) defines an analytic function in 
the half plane (\ref{U243}).  The function on the left of (\ref{U245}) is obviously analytic in the half plane (\ref{U243}), and the 
function on the right of (\ref{U245}) is analytic in the half plane (\ref{U243}) by our previous remarks.  It follows from (\ref{U375}) that 
these functions are equal on the positive real axis.  Hence they are equal in the half plane (\ref{U243}) by analytic continuation.
\end{proof}

We now consider a more general situation.  Let $L_1, L_2, \dots , L_J$, be positive integers, and for each index $j = 1, 2, \dots , J$, let
\begin{equation*}\label{U193}
\delta_{L_j} : \R^{L_j} \rightarrow [0, \infty)
\end{equation*}
be a distance function.  Write
\begin{equation}\label{U200}
\fS_{L_j} = \big\{\bx_j \in \R^{L_j} : \delta_{L_j}\bigl(\bx_j\bigr) \le 1\big\}
\end{equation}
for the closed star-body in $\R^{L_j}$ determined by $\delta_{L_j}$, and assume that 
\begin{equation}\label{U207}
\Vol_{L_j}\bigl(\fS_{L_j}\bigr) < \infty,\quad\text{for each $j = 1, 2, \dots , J$}.
\end{equation}
If $0 < \xi < \infty$ then as in (\ref{U16}) we have
\begin{equation*}\label{U214}
\Vol_{L_j}\big\{\bwy \in \R^{L_j} : \delta_{L_j}\bigl(\bwy\bigr) \le \xi\big\} = \Vol_{L_j}\bigl(\fS_{L_j}\bigr) \xi^{L_j}.
\end{equation*}  
Next we define
\begin{equation}\label{U221}
N = L_1 + L_2 + \cdots + L_J,
\end{equation}
and we write $\R^N$ as
\begin{equation}\label{U228}
\R^N = \big\{\bigl(\bx_1, \bx_2, \dots , \bx_J\bigr) : \text{$\bx_j \in \R^{L_j}$ for $j = 1, 2, \dots , J$}\big\}.
\end{equation}
We use the collection of distance functions $\delta_{L_j}$ for $j = 1, 2, \dots , J$, to define a distance function
\begin{equation*}\label{U235}
\Delta_N : \R^N \rightarrow [0, \infty)
\end{equation*}
by
\begin{equation}\label{U242}
\Delta_N\big(\bx_1, \bx_2, \dots , \bx_J\bigr) = \delta_{L_1}\bigl(\bx_1\bigr) + \delta_{L_2}\bigl(\bx_2\bigr) + \cdots + \delta_{L_J}\bigl(\bx_J\bigr).
\end{equation}
We write
\begin{equation*}\label{U244}
\fS_N = \big\{\bigl(\bx_1, \bx_2, \dots , \bx_J\bigr) \in \R^N : \Delta_N\big(\bx_1, \bx_2, \dots , \bx_J\bigr) \le 1\big\}
\end{equation*}
for the closed star-body in $\R^N$ determined by $\Delta_N$.  As
\begin{equation*}\label{U246}
\begin{split}
\fS_N &\subseteq \big\{\bigl(\bx_1, \bx_2, \dots , \bx_J\bigr) \in \R^N : \text{$\delta_{L_j}\bigl(\bx_j\bigr) \le 1$ for $j = 1, 2, \dots , J$}\big\}\\
          &= \prod_{j = 1}^J \big\{\bx_j \in \R^{L_j} : \delta_{L_j}\bigl(\bx_j\bigr) \le 1\big\}
          = \fS_{L_1} \times \fS_{L_2} \times \cdots \times \fS_{L_J},
\end{split}
\end{equation*}
it follows from (\ref{U207}) that $\Vol_N\bigl(\fS_N\bigr)$ is finite.

\begin{theorem}\label{thmU1}  Let $L_1, L_2, \dots , L_J$, be positive integers, and let
\begin{equation*}\label{U249}
N = L_1 + L_2 + \cdots + L_J.
\end{equation*}
For each integer $j = 1, 2, \dots , J$, assume that $\delta_{L_j}$ is a distance function on $\R^{L_j}$ such that the star-body $\fS_{L_j}$ 
defined by {\rm (\ref{U200})} has finite volume.  Let $\Delta_N$ be the distance function on $\R^N$ defined by {\rm (\ref{U242})}. 
Then the volume of the star-body
\begin{equation*}\label{U256}
\fS_N = \big\{\bigl(\bx_1, \bx_2, \dots , \bx_J\bigr) \in \R^N : \Delta_N\big(\bx_1, \bx_2, \dots , \bx_J\bigr) \le 1\big\}
\end{equation*}
is finite and satisfies the identity
\begin{equation}\label{U263}
N! \Vol_N\bigl(\fS_N\bigr) = \prod_{j = 1}^J \bigl(L_j ! \Vol_{L_j}\bigl(\fS_{L_j}\bigr)\bigr).
\end{equation}
\end{theorem}

\begin{proof}  We have already noted that $\fS_N$ has finite volume in $\R^N$.  Therefore it follows from Lemma \ref{lemU1} that
\begin{equation}\label{U271}
N! \Vol_N(\fS) s^{-N} = \int_{\R^N} \exp\bigl(-s \Delta(\bwy)\bigr)\ \dby
\end{equation}
at each point $s = \sigma + it$ such that $0 < \sigma$.  Using (\ref{U242}) and Lemma \ref{lemU1} again we find that  
\begin{equation}\label{U407}
\begin{split}
 \int_{\R^N} \exp\bigl(- s \Delta_N(\bwy)\bigr)\ \dby &= \prod_{j = 1}^J \int_{\R^{L_j}} \exp\bigl(- s \delta_{L_j}(\bx_j)\bigr)\ \dbx_j\\
 	&= \prod_{j = 1}^J \bigl(L_j! \Vol_{L_j}\bigl(\fS_{L_j}\bigr) s^{-L_j}\bigr)\\
	&= \biggl(\prod_{j = 1}^J L_j! \Vol_{L_j}\bigl(\fS_{L_j}\bigr)\biggr) s^{-N}.
\end{split}	
\end{equation} 
The identity (\ref{U263}) follows from (\ref{U271}) and (\ref{U407}).
\end{proof}

In our application of (\ref{U263}) in section \ref{heights} the distance functions that occur are norms.  Thus the star-bodies that occur 
are compact, convex, symmetric subsets of $\R^N$.

%%%%%%%%%%%%%%%%%%%%%%%%%%%%%%%%%%%%%%%%%%%%%%%%%%%%%%%%%%%%%%%%%%%%
\section{Heights, relative units, and the relative regulator}\label{heights}

At each place $v$ of $k$ we write $k_v$ for the completion of $k$ at $v$.  We work with two distinct absolute values $\|\ \|_v$ and 
$|\ |_v$ from each place $v$.  These absolute values are related by
\begin{equation*}\label{B4}
\|\ \|_v^{d_v/d} = |\ |_v,
\end{equation*}
where $d_v = [k_v : \Q_v]$ is the local degree at $v$, and $d = [k : \Q]$ is the global degree.  If $v | \infty$ then the restriction of $\|\ \|_v$ to 
$\Q$ is the usual archimedean absolute value on $\Q$, and if $v | p$ then the restriction of $\|\ \|_v$ to $\Q$ is the usual $p$-adic absolute 
value on $\Q$.  The absolute logarithmic Weil height is the map
\begin{equation}\label{B6}
h : k^{\times} \rightarrow [0, \infty)
\end{equation}
defined at each algebraic number $\alpha \not= 0$ in $k$ by the sum
\begin{equation}\label{B12}
h(\alpha) = \sum_v \log^+ |\alpha|_v = \hh \sum_v \bigl|\log |\alpha|_v \bigr|.
\end{equation}
In both sums there are only finitely many nonzero terms, and the equality on the right of (\ref{B12}) follows from the product formula.  It can be 
shown that the value of $h(\alpha)$ does not depend on the field $k$ that contains $\alpha$.  Hence the Weil height may be regarded as a map
\begin{equation*}\label{B14}
h : \oq \rightarrow [0, \infty),\quad\text{or as a map}\quad h : \oQt \rightarrow [0, \infty).
\end{equation*}
Further information about the height is contained in \cite[section 1.5]{bombieri2006}.

Let the number fields $k$ and $l$ satisfy the hypotheses (\ref{short277}) and (\ref{short298}).  We recall (see 
\cite[equation (3.1)]{akhtari2016}) that at each place $v$ of $k$ and for each element $\alpha$ in $l^{\times}$, we have
\begin{equation}\label{extra271}
[l : k] \sum_{w | v} \log |\alpha|_w = \log |\Norm_{l/k}(\alpha)|_v,
\end{equation}
where the sum on the left of (\ref{extra271}) is over the set of places $w$ of $l$ such that $w | v$.  In particular, it follows from 
the definition (\ref{extra259}) that $\alpha$ in $F_l$ belongs to the kernel $E_{l/k}$ if and only if
\begin{equation}\label{extra263}
\sum_{w | v} [l_w : \Q_w] \log \|\alpha\|_w = [k_v : \Q_v] \log \|\norm_{l/k}(\alpha)\|_v = 0
\end{equation}
at each archimedean place $v$ of $k$.  If $\alpha$ belongs to $F_l$ then the product formula implies that
\begin{equation}\label{extra273}
\sum_{w | \infty} [l_w : \Q_w] \log \|\alpha\|_w = 0,
\end{equation}
where the sum on the left of (\ref{extra273}) is over the set of all archimedean places $w$ of $l$.  If $\alpha$ is a relative unit in $E_{l/k}$ 
then there are subsums of (\ref{extra273}) that equal zero and are given by (\ref{extra263}) for each archimedean place $v$ of $k$.  

At each place $v$ of $k$ we define
\begin{equation*}\label{extra275}
W_v(l/k) = \{w : \text{$w$ is a place of $l$, and $w|v$}\}.
\end{equation*}
We recall that $r(k) + 1$ is the number of archimedean places of $k$, and $r(l) + 1$ is the number of archimedean places of $l$.
At each archimedean place $v$ of $k$ we select a place $\widehat{w}_v$ of $l$ such that $\widehat{w}_v | v$.  Then it follows that the set
of archimedean places
\begin{equation}\label{extra280}
S(l/k) = \bigcup_{v|\infty} \bigl(W_v(l/k) \setminus \{\widehat{w}_v\}\bigr)
\end{equation}
has cardinality
\begin{equation}\label{extra285}
\begin{split}
|S(l/k)| &= \sum_{v | \infty} \Bigl(\bigl|W_v(l/k)\bigr| - |\{\widehat{w}_v\}|\Bigr)\\
            &= \bigl(r(l) + 1\bigr) - \bigl(r(k) + 1\bigr)\\
	     &= r(l/k).
\end{split}
\end{equation}
The disjoint union on the right of (\ref{extra280}) provides a partition of the set $S(l/k)$ into disjoint subsets indexed by the 
collection of all archimedean places $v$ of $k$.  For each archimedean place $v$ of $k$ we define
\begin{equation*}\label{extra286}
T_v = W_v(l/k) \setminus \{\widehat{w}_v\},
\end{equation*}  
so that (\ref{extra280}) can be written as the disjoint union
\begin{equation}\label{extra287}
S(l/k) = \bigcup_{v | \infty} T_v.
\end{equation}

Let $\R^{S(l/k)}$ be the real vector space of functions $\bx : S(l/k) \rightarrow \R$.  And for each archimedean place $v$ of $k$ let
$\R^{T_v}$ be the subspace of functions $\bx$ in $\R^{S(l/k)}$ that have support contained in the subset $T_v$. 
The Schinzel norm (see \cite[equation (4.3)]{akhtari2016}) is the map
\begin{equation*}\label{V7}
\delta : \R^{S(l/k)} \rightarrow [0, \infty)
\end{equation*}
defined by
\begin{equation}\label{V14}
\delta(\bx) = \hh \biggl|\sum_{w \in S(l/k)} x_w\biggr| + \hh \sum_{w \in S(l/k)} |x_w|.
\end{equation}
We write
\begin{equation}\label{extra288}
\pi_{T_v} : \R^{S(l/k)} \rightarrow \R^{T_v}
\end{equation}
for the linear projection of $\R^{S(l/k)}$ onto $\R^{T_v}$.  Then for each archimedean place $v$ of $k$ the composition
\begin{equation}\label{V28}
\delta \circ \pi_{T_v} : \R^{S(l/k)} \rightarrow [0, \infty)
\end{equation}
is given by
\begin{equation}\label{V35}
\delta\bigl(\pi_{T_v}(\bx)\bigr) = \hh \biggl|\sum_{w \in T_v} x_w\biggr| + \hh \sum_{w \in T_v} |x_w|.
\end{equation}
We use the projections (\ref{extra288}) to define the {\it generalized Schinzel norm} associated to the partition (\ref{extra287}).  More
precisely, we define the generalized Schinzel norm
\begin{equation*}\label{extra289}
\Delta : \R^{S(l/k)} \rightarrow [0, \infty),
\end{equation*}
to be the map
\begin{equation}\label{extra290}
\begin{split}
\Delta(\bx) = \sum_{v | \infty} \delta\bigl(\pi_{T_v}(\bx)\bigr) 
	= \hh \sum_{v | \infty} \biggl|\sum_{w \in T_v} x_w\biggr| + \hh \sum_{v | \infty}  \sum_{w \in T_v} |x_w|.
\end{split}
\end{equation}

For each place $v$ of $k$ the restriction of (\ref{V28}) to the subspace $\R^{T_v} \subseteq \R^{S(l/k)}$ is a norm on $\R^{T_v}$.  
The unit ball associated to this norm is the compact, convex, symmetric subset
\begin{equation*}\label{V42}
\fS_v = \big\{\bx \in \R^{T_v} : \delta\bigl(\pi_{T_v}(\bx)\bigr) \le 1\big\}.
\end{equation*}
It was shown in \cite[Lemma 4.1]{akhtari2016} that
\begin{equation}\label{V49}
\Vol_{|T_v|}\bigl(\fS_v\bigr) = \frac{(2 |T_v|)!}{(|T_v|!)^3},\quad\text{or}\quad |T_v|! \Vol_{|T_v|}\bigl(\fS_v\bigr) = \binom{2 |T_v|}{|T_v|},
\end{equation}
where the expression on the right of (\ref{V49}) is a middle binomial coefficient.  Let
\begin{equation}\label{V56}
\fS_{l/k} = \big\{\bx \in \R^{S(l/k)} : \Delta(\bx) \le 1\big\}
\end{equation}
be the closed unit ball defined using the generalized Schinzel norm (\ref{extra290}).  As $\Delta$ is a norm on $\R^{S(l/k)}$, it follows that
$\fS_{l/k}$ is also a compact, convex, symmetric subset.  And it follows from Theorem \ref{thmU1} and (\ref{extra285}) that 
\begin{equation}\label{V63}
r(l/k)! \Vol_{r(l/k)}\bigl(\fS_{l/k}\bigr) = \prod_{v | \infty} \Bigl(|T_v|! \Vol_{|T_v|}\bigl(\fS_v\bigr)\Bigr)
	= \prod_{v | \infty} \binom{2 |T_v|}{|T_v|}.
\end{equation}
The product on the right of (\ref{V63}) is not easy to express using the invariants $r(k)$ and $r(l)$.  However, we have the simple inequality
\begin{equation}\label{V70}
2^{r(l/k)} = \prod_{v | \infty} 2^{|T_v|} \le \prod_{v | \infty} \binom{2 |T_v|}{|T_v|} = r(l/k)! \Vol_{r(l/k)}\bigl(\fS_{l/k}\bigr).
\end{equation}
We use (\ref{V70}) in our proof of Theorem \ref{thmintro1}, but more elaborate lower bounds involving $r(k)$ and $r(l)$ can also be proved.

Let $\eta_1, \eta_2, \dots , \eta_{r(l/k)}$ be a collection of multiplicatively independent relative units that form a basis for the
subgroup $E_{l/k} \subseteq F_l$ of relative units.  Let  $M_{l/k}$ be the $r(l/k) \times r(l/k)$ real matrix
\begin{equation}\label{extra294}
M_{l/k} = \bigl([l_w : \Q_w] \log \|\eta_j\|_w\bigr),
\end{equation}
where $w \in S(l/k)$ indexes rows and $j = 1, 2, \dots , r(l/k)$ indexes columns.  As in \cite[equation (3.5)]{akhtari2016} and 
\cite[equation (1.3)]{costa1991}, the {\it relative regulator} of $l/k$ is the positive number
\begin{equation}\label{extra301}
\Reg\bigl(E_{l/k}\bigr) = \bigl|\det M_{l/k}\bigr|.
\end{equation}
It follows, as in the proof of \cite[Theorem 1]{costa1991} (see also \cite{costa1993}), that the absolute value of the determinant on the right 
of (\ref{extra301}) does not depend on the choice of places $\widehat{w}_v$ that are removed from each subset $W_v(l/k)$ for each 
archimedean place $v$ of $k$.  Alternatively, the (ordinary) regulators $\Reg(k)$ and $\Reg(l)$, and the relative regulator 
$\Reg\bigl(E_{l/k}\bigr)$, are related (see \cite[Theorem 1]{costa1991}) by the identity
\begin{equation}\label{extra302}
[F_k : I_{l/k}] \Reg(k) \Reg\bigl(E_{l/k}\bigr) = \Reg(l),
\end{equation}
where $I_{l/k}$ is the subgroup defined in (\ref{extra277}).

The following result is a refinement of \cite[Lemma 5.1]{akhtari2016} applicable to relative units.

\begin{lemma}\label{lemheight1}  Let $\alpha$ belong to the group $F_l$, and let 
\begin{equation}\label{short300}
\bx(\alpha) = \bigl([l_w : \Q_w] \log \|\alpha\|_w\bigr)
\end{equation}
be the image of $\alpha$ in the real vector space $\R^{S(l/k)}$, where the rows of the column vector $\bx(\alpha)$ are indexed by places 
$w$ in $S(l/k)$.  If $\alpha$ belongs to the subgroup $E_{l/k}$ of relative units, then
\begin{equation}\label{short305}
\Delta\bigl(\bx(\alpha)\bigr) = [l : \Q] h(\alpha).
\end{equation}
\end{lemma}

\begin{proof}  We recall from (\ref{extra263}) that for each archimedean place $v$ of $k$ we have
\begin{equation}\label{short310}
[l_{\widehat{w}_v} : \Q_{\widehat{w}_v}] \log \|\alpha\|_{\widehat{w}_v} 
	+ \sum_{\substack{w | v\\w \not= \widehat{w}_v}} [l_w : \Q_w] \log \|\alpha\|_w = 0.
\end{equation}
Using (\ref{V35}) and (\ref{short310}) we find that
\begin{equation}\label{short314}
\begin{split}
2 \delta\bigl(\pi_{T_v}(\bx(\alpha)\bigr) &= \biggl|\sum_{\substack{w | v\\w \not= \widehat{w}_v}} [l_w : \Q_w] \log \|\alpha\|_w\biggr|
				+ \sum_{\substack{w | v\\w \not= \widehat{w}_v}} [l_w : \Q_w] \bigl|\log \|\alpha\|_w\bigr|\\
				&= [l_{\widehat{w}_v} : \Q_{\widehat{w}_v}] \bigl|\log \|\alpha\|_{\widehat{w}_v}\bigr|
				+ \sum_{\substack{w | v\\w \not= \widehat{w}_v}} [l_w : \Q_w] \bigl|\log \|\alpha\|_w\bigr|\\
                                  &= \sum_{w|v} [l_w : \Q_w] \bigl|\log \|\alpha\|_w\bigr|\\
                                  &= [l : \Q] \sum_{w | v} \bigl|\log |\alpha|_w\bigr|.
\end{split}
\end{equation}
Next we combine (\ref{extra290}) and (\ref{short314}) to get
\begin{equation}\label{short320}
\begin{split}
\Delta\bigl(\bx(\alpha)\bigr) &= \hh [l : \Q] \sum_{v | \infty} \sum_{w | v} \bigl|\log |\alpha|_w\bigr|\\
					&= \hh [l : \Q] \sum_{w | \infty} \bigl|\log |\alpha|_w\bigr|\\
					&= [l : \Q] h(\alpha).
\end{split}
\end{equation}
This proves (\ref{short305}).
\end{proof}

%%%%%%%%%%%%%%%%%%%%%%%%%%%%%%%%%%%%%%%%%%%%%%%%%%%%%%%%%%%%%%%%%%
\section{The proof of Theorem \ref{thmintro1} and Corollary \ref{corintro1}}\label{proofs}

As before, let $\eta_1, \eta_2, \dots , \eta_{r(l/k)}$ be a basis for the subgroup $E_{l/k} \subseteq F_l$ of relative units. 
Then let $\vep_1, \vep_2, \cdots , \vep_{r(l/k)}$ be a basis for the subgroup $\fE \subseteq E_{l/k}$.  
It follows that there exists an $r(l/k) \times r(l/k)$ nonsingular matrix $C = (c_{i j})$ with entries in $\Z$ such that
\begin{equation}\label{short311}
\log \|\vep_j\|_w = \sum_{i = 1}^{r(l/k)} c_{i j} \log \|\eta_i\|_w
\end{equation}
at each archimedean place $w$ of $l$.  The system of equations (\ref{short311}) can be written as the single matrix equation
\begin{equation}\label{short315}
\bigl([l_w : \Q_w]\log \|\vep_j\|_w\bigr) = \bigl([l_w : \Q_w]\log \|\eta_j\|_w\bigr) C,
\end{equation}
where $w$ is an archimedean place of $l$ and $w$ indexes the rows of the matrices which are not $C$ on both sides of (\ref{short315}).

At each archimedean place $v$ of $k$ we remove the row indexed by $\widehat{w}_v$ in the matrix on the left of (\ref{short315})
so as to obtain an $r(l/k) \times r(l/k)$ submatix
\begin{equation*}\label{short321}
L(\fE) = \bigl([l_w : \Q_w]\log \|\vep_j\|_w\bigr).
\end{equation*}
We note that the rows of $L(\fE)$ are indexed by the places $w$ in the set $S(l/k)$.
Removing the same rows in the product on the right of (\ref{short315}) leads to the matrix identity
\begin{equation}\label{short325}
L(\fE) = M_{l/k} C,
\end{equation}
where $M_{l/k}$ was defined in (\ref{extra294}).  We conclude from (\ref{short325}) that
\begin{equation}\label{short330}
|\det L(\fE)| = \Reg(E_{l/k}) [E_{l/k} : \fE].
\end{equation}

We use the $r(l/k) \times r(l/k)$ matrix $L(\fE)$ to define a lattice
\begin{equation*}\label{short335}
\eL = \big\{L(\fE) \bxi : \bxi \in \Z^{S(l/k)}\big\} \subseteq \R^{S(l/k)}.
\end{equation*}
The determinant of $\eL$ is given by (\ref{short330}).  Let $\Delta$ be the norm defined by (\ref{extra290}), and let
\begin{equation*}\label{short340}
0 < \lambda_1 \le \lambda_2 \le \cdots \le \lambda_{r(l/k)} < \infty
\end{equation*}
be the successive minima associated to the lattice $\eL$ and the compact, convex symmetric set
\begin{equation*}\label{short345}
\fS_{l/k} = \big\{\bx \in \R^{S(l/k)} : \Delta(\bx) \le 1\big\}.
\end{equation*}
Then there exist multiplicatively independent points $\psi_1, \psi_2, \dots , \psi_{r(l/k)}$ in $\fE$ such that
\begin{equation*}\label{short355}
\Delta\bigl(\bx(\psi_j)\bigr) = \lambda_j,\quad\text{for each $j = 1, 2, \dots , r(l/k)$}.
\end{equation*}
 As in (\ref{short300}) we have written
\begin{equation}\label{short360}
\bx(\psi_j) = \bigl([l_w : \Q_w] \log \|\psi_j\|_w\bigr),
\end{equation}
where the rows of the column vector on the right of (\ref{short360}) are indexed by places $w$ in $S(l/k)$.
From Minkowski's theorem (see \cite[Chapter VIII.4.3, Theorem V]{cassels1971}) on successive minima we get the inequality
\begin{equation}\label{short350}
\Vol_{S(l/k)}\bigl(\fS_{l/k}\bigr) \lambda_1 \lambda_2 \cdots \lambda_{r(l/k)} \le 2^{r(l/k)}\Reg(E_{l/k}) [E_{l/k} : \fE].
\end{equation}
And from Lemma \ref{lemheight1} we find that
\begin{equation}\label{short365}
\Delta\bigl(\bx(\psi_j)\bigr) = \lambda_j = [l : \Q] h(\psi_j)
\end{equation}
for each $j = 1, 2, \dots , r(l/k)$.  Combining (\ref{short350}) and (\ref{short365}) leads to the inequality
\begin{equation}\label{short371}
\prod_{j = 1}^{r(l/k)} \bigl([l : \Q] h(\psi_j)\bigr) 
	\le \frac{2^{r(l/k)} \Reg\bigl(E_{l/k}\bigr) [E_{l/k} : \fE]}{\Vol_{S(l/k)}\bigl(\fS_{l/k}\bigr)}.
\end{equation}
Finally, it follows from the inequality (\ref{V70}) that
\begin{equation}\label{short373}
2^{r(l/k)} \Bigl(\Vol_{S(l/k)}\bigl(\fS_{l/k}\bigr)\Bigr)^{-1} \le r(l/k)!
\end{equation}
Then (\ref{short371}) and (\ref{short373}) lead to the inequality (\ref{short322}) in the statement of the theorem.
\smallskip

We now prove Corollary \ref{corintro1}.  Assume that the elements in the set (\ref{short450}) are multiplicatively dependent.  Then there 
exist lattice points $\bm$ in $\Z^{r(k)}$ and $\bn$ in $\Z^{r(l/k)}$, such that $\bm$ and $\bn$ are not both $\bo$ and
\begin{equation}\label{short460}
\beta_1^{m_1} \beta_2^{m_2} \cdots \beta_{r(k)}^{m_{r(k)}} \psi_1^{n_1} \psi_2^{n_2} \cdots \psi_{r(l/k)}^{n_{r(l/k)}} = 1.
\end{equation}
If $\bm = \bo$ then it follows from the independence of $\psi_1, \psi_2, \dots , \psi_{r(l/k)}$ that $\bn = \bo$.  Similarly,
if $\bn = \bo$ then it follows that $\bm = \bo$.  Thus we may assume that both $\bm \not= \bo$ and $\bn \not= \bo$.  

Let
\begin{equation*}\label{short465}
\alpha = \beta_1^{m_1} \beta_2^{m_2} \cdots \beta_{r(k)}^{m_{r(k)}},\quad\text{and}\quad
		\gamma = \psi_1^{n_1} \psi_2^{n_2} \cdots \psi_{r(l/k)}^{n_{r(l/k)}}.
\end{equation*}
Then $\alpha \not= 1$ belongs to $I_{l/k}$ and $\gamma \not= 1$ belongs to $E_{l/k}$, and it follows from (\ref{short460}) that
$\alpha = \gamma^{-1}$.  Thus $\alpha \not= 1$ belongs to both $I_{l/k}$ and $E_{l/k}$, which is impossible by (\ref{extra270}).  We
have verified that the elements of the set (\ref{short450}) are multiplicatively independent.    
The inequality (\ref{short455}) follows now by using (\ref{short440}), (\ref{short445}), and (\ref{extra302}).

%%%%%%%%%%%%%%%%%%%%%%%%%%%%%%%%%%%%%%%%%%%%%%%%%%%%%%%
\section*{Acknowledgements} The authors acknowledge the support from \emph{Max Planck Institute for Mathematics in Bonn}, where this project was initiated.  The authors are grateful to the anonymous referee for insightful comments. Shabnam Akhtari's research has been in parts supported by  \emph{the Simons Foundation Collaboration Grants}, Award Number 635880, and by  \emph{the National Science Foundation} Award DMS-2001281.
 
%%%%%%%%%%%%%%%%%%%%%%%%%%%%%%%%%%%%%%%%%%%%%%%%%%%%%%%%%%%%%%%

%\today
\end{document}